\documentclass[11pt]{article}
\usepackage{amsmath,amssymb} 

\begin{document}
\title{Immersed surfaces and Seifert fibered \\
 surgery on Montesinos knots} 
\author{Ying-Qing Wu} 
\date{}
\maketitle

\footnotetext[1]{ Mathematics subject classification:  {\em Primary 
57N10.}}

\footnotetext[2]{ Keywords and phrases: Immersed surfaces, Dehn
  surgery, Seifert fibered manifolds, Montesinos knots}

\begin{abstract}
  We will use immersed surfaces to study Seifert fibered surgery on
  Montesinos knots, and show that if $\frac 1{q_1-1} + \frac 1{q_2-1}
  + \frac 1{q_3-1} \leq 1$ then a Montesinos knot $K(\frac{p_1}{q_1},
  \frac{p_2}{q_2}, \frac{p_3}{q_3})$ admits no atoroidal Seifert
  fibered surgery.
\end{abstract}

\newcommand{\proof}{\noindent {\bf Proof.} }
\newcommand{\qed}{\quad $\Box$}
\newtheorem{thm}{Theorem}[section]
\newtheorem{prop}[thm]{Proposition} 
\newtheorem{lemma}[thm]{Lemma} 
\newtheorem{cor}[thm]{Corollary} 
\newtheorem{defn}[thm]{Definition} 
\newtheorem{notation}[thm]{Notation} 
\newtheorem{qtn}[thm]{Question} 
\newtheorem{example}[thm]{Example} 
\newtheorem{remark}[thm]{Remark} 
\newtheorem{conj}[thm]{Conjecture} 
\newtheorem{prob}[thm]{Problem} 
\newtheorem{rem}[thm]{Remark} 

\newcommand{\bdd}{\partial}
\newcommand{\Int}{{\rm Int}}
\renewcommand{\a}{\alpha}
\renewcommand{\b}{\beta}
\newcommand{\e}{\epsilon}

\input epsf.tex

\section{Introduction}

Exceptional Dehn surgeries on arborescent knots have been studied
extensively.  They have been classified for arborescent knots of
length at least 4 [Wu2], as well as for all 2-bridge knots [BW].
There is no reducible surgery on hyperbolic arborescent knots [Wu1],
and toroidal surgeries on length 3 Montesinos knots have also been
classified [Wu3].  Therefore atoroidal Seifert fibered surgeries on
Montesinos knots of length 3 are the only ones on arborescent knots
that have not been determined.

Atoroidal Seifert fibered surgery is much more difficult to deal with
than other types of exceptional surgeries.  For example, the minimum
upper bounds for distances between two types of exceptional Dehn
fillings on hyperbolic manifolds have all been determined when those
are not atoroidal Seifert fibered, but no such bound is known when one
of them is.  See [GW] and the references there for works in that
direction.  The major difficulty to deal with atoroidal Seifert
fibered surgery is that there is no embedded essential small surfaces
(sphere, disk, annulus or torus) in such manifolds and therefore one
cannot use those traditional combinatorial methods on intersection
graphs for such surgery problems.  Those with infinite fundamental
group do contain immersed essential tori, but there had not been much
success using them as tools in solving Dehn surgery problems.  In this
paper, however, we {\em will\/} use those immersed surfaces as a major
tool.  We will study the intersection of an immersed surface $F$ with
the tangle decomposition surfaces and the tangle spaces, and show that
no such surface could exist when the tangles are not too simple.

A length 3 Montesinos knot $K(\frac{p_1}{q_1}, \frac{p_2}{q_2},
\frac{p_3}{q_3})$ is the cyclic union of three rational tangles
$T(\frac{p_i}{q_i})$, where $q_i \geq 2$ and $p_i, q_i$ are coprime.
The integer part of the $p_i/q_i$ can be shifted around, so we will
always assume that $2|p_i| \leq q_i$ for $i=1,2$.  Given a knot $K$ in
$S^3$, we use $K(r)$ to denote the manifold obtained by Dehn surgery
on $K$ along a slope $r$ on $\bdd N(K)$, where $N(K)$ is a tubular
neighborhood of $K$.  The following is our main theorem.

%1.1
\begin{thm} Suppose $K = K(\frac{p_1}{q_1}, \frac{p_2}{q_2},
  \frac{p_3}{q_3})$ is a hyperbolic Montesinos knot of length 3.  If
  $\frac 1{q_1-1} + \frac 1{q_2-1} + \frac 1{q_3-1} \leq 1$ then $K$
  admits no atoroidal Seifert fibered surgery.
\end{thm}

Recently Ichihara and Jong [IJ2] showed that the only toroidal Seifert
fibered surgery on Montesinos knots is the $0$ surgery on the trefoil
knot, hence the theorem is still true with the word ``atoroidal''
deleted.  By [Wu1] there is no reducible surgery on hyperbolic
Montesinos knots, so the following result follows immediately from
Theorem 1.1 and the classification of toroidal surgeries on these
knots [Wu3, Theorem 1.1].

%1.2
\begin{cor} Suppose $K = K(\frac{p_1}{q_1}, \frac{p_2}{q_2},
  \frac{p_3}{q_3})$ is a hyperbolic Montesinos knot of length 3.  If
  $\frac 1{q_1-1} + \frac 1{q_2-1} + \frac 1{q_3-1} \leq 1$ and a Dehn
  surgery $K(r)$ is nonhyperbolic, then $|p_i| = 1$, $r$ is the
  pretzel slope, and $K(r)$ is toroidal.
\end{cor}

We may assume $2 \leq q_1 \leq q_2 \leq q_3$.  Corollary 1.2 and [BW,
Wu2] provide a classification of exceptional Dehn surgeries on all
arborescent knots except those $K(\frac{p_1}{q_1}, \frac{p_2}{q_2},
\frac{p_3}{q_3})$ with $q_1 = 2$, or $(q_1, q_2) = (3,3)$, or $(q_1,
q_2, q_3) = (3,4,5)$.  Further restrictions on $p_i$ and the surgery
slopes for these cases will be given in Theorems 8.2.

Theorem 1.1 is known in the special case that $K(r)$ has finite
fundamental group.  The classification of finite surgeries on
Montesinos knots has been completed by Ichihara and Jong [IJ1].  It
used the results of Delman [De] on essential laminations, Mattman's
result [Ma] on surgery on pretzel knots, and Ni's result [Ni] on
Heegaard Floer homology and fibered knots.  See also [Wa] and [FIKMS].
While our major goal is to prove Theorem 1.1 for infinite Seifert
fibered surgeries, we will also provide an independent proof of it for
the finite surgery case.  We extend the thin position idea of Gabai
[Ga] and use an immersed thin sphere to replace the immersed essential
torus in the case that $K(r)$ has infinite fundamental group.  This
works fine in our setting and sometimes the proof is simpler than for
infinite surgery case since the surface is now a sphere instead of a
torus.

The paper is organized as follows.  Section 2 is to set up some
notations and conventions, and introduce some basic lemmas.  Section 3
discusses immersed essential disks in tangle spaces.  It will be shown
that any such disk must intersect the axis of the tangle at some
minimal number of points, and the disks are embedded and standard in
certain cases.  In Section 4 we define immersed surfaces that are in
essential position with respect to an essential embedded surface and
prove its existence in manifolds with finite fundamental groups.
Section 5 defines elementary surfaces and shows that if the surface
$F$ coming from an immersed $\pi_1$-injective torus is elementary then
the surgered manifold is either toroidal or the connected sum of two
lens spaces.  This is crucial to deal with the fact that some of the
knots excluded in Theorems 1.1 and 8.2 admit toroidal surgeries and
hence the surgered manifold does contain $\pi_1$-injective tori.  In
Section 6 we define intersection graphs and prove some basic
properties of such graphs.  Section 7 defines angled Euler numbers and
show that it is additive.  Section 8 completes the proof of the main
theorems.

\section{Preliminaries}

In this paper we will consider both embedded surfaces and non-embedded
surfaces in 3-manifolds.  We always assume that surfaces and curves
intersect transversely.  Unless otherwise stated, surfaces $F$ in a
3-manifold $M$ are assumed to have boundaries on the boundary of $M$,
and a homotopy of $F$ refers to a relative homotopy of the pair $(F,
\bdd F)$ in $(M, \bdd M)$, i.e.\ a homotopy $f_t: F\to M$ such that
$f_t(\bdd F) \subset \bdd M$ for all $t$.  Similarly for arcs on
surfaces.

An arc $\a$ on a surface $F$ is trivial if it is rel $\bdd \a$
homotopic to an arc on $\bdd F$.  If $\a$ is closed then it is trivial
if and only if it is null homotopic on $F$.  A (possibly non-embedded)
disk $D$ in a 3-manifold $M$ is {\it nontrivial\/} if it is not rel
$\bdd D$ homotopic to a disk on $\bdd M$.  If $M$ is irreducible (in
particular if $M$ is the tangle space $E(t)$ below), then $D$ in $M$
is nontrivial if and only if $\bdd D$ is a nontrivial curve on $\bdd
M$.  A curve or disk is {\it essential\/} if it is nontrivial.  Two
(possibly immersed) curves $C_1, C_2$ on a surface $F$ {\it intersect
  minimally\/} if there are no subarc $a_i \subset C_i$ such that
$\bdd a_1 = \bdd a_2$ and the loop $a_1 \cup a_2$ is null homotopic on
$F$.  When $C_1, C_2$ are embedded, this is equivalent to say that
$C_1 \cup C_2$ contains no bigons on $F$, i.e.\ there is no arcs
$a_i\subset C_i$ with $\bdd a_1 = \bdd a_2$ such that $a_1 \cup a_2$
bounds a disk on $F$ with interior disjoint from $C_1 \cup C_2$.

In this paper a {\it tangle\/} is a triple $T = (B, t, m)$, where $B$
is a fixed 3-ball, $t = t_1 \cup t_2$ is a pair of arcs properly
embedded in $B$, and $m$ is a simple loop on $\bdd B$, cutting $\bdd
B$ into two disks, called the left disk and the right disk, each
containing two points of $\bdd t$.  The curve $m$ is called the {\it
  axis\/} of the tangle.  Two tangles $(B, t, m)$ and $(B', t', m')$
are {\it equivalent\/} if they are homeomorphic as a triple.  They are
{\it strongly equivalent\/} if $B=B'$ and the homeomorphism is the
identity on $\bdd B$.

Denote by $N(t)$ a regular neighborhood of $t$, and by $E(t) = B - \Int
N(t)$ the exterior of $t$, which will also be called the {\it tangle
  space\/} of $(B, t, m)$.  Denote by $A(t)$ the two annuli $A(t) =
\bdd N(t) \cap E(t) = A_1(t) \cup A_2(t)$ on $\bdd E(t)$, by $P(t)$
the $4$-punctured sphere $\bdd B \cap E(t)$, and by $P_1(t) \cup
P_2(t)$ the two twice punctured disks obtained by cutting $P(t)$ along
$m$.  If $C$ is a properly embedded $n-1$ manifold in an $n$-manifold
$F$, denote by $F|C$ the manifold obtained by cutting $F$ along $C$.

%2.1
\begin{defn} {\rm A homotopy or isotopy $h_x$ of $E(t)$ or $\bdd E(t)$
    is {\em $P$-preserving\/} if $h_x$ maps each of the set $A(t),
    P_1(t), P_2(t), m$ to itself during the homotopy.  Similarly, if
    $C$ is a curve or surface in $E(t)$ then a $P$-preserving homotopy
    or isotopy $h_x$ of $C$ is such that $C \cap Y$ is mapped to $Y$
    for $Y = A(t), P_1(t), P_2(t), m$ and all $x \in [0,1]$.}
\end{defn}

A tangle $T=(B,t,m)$ is a $p/q$ {\it rational tangle\/}, denoted by
$T(p/q)$, if $B$ is isotopic to a pillowcase with the four points of
$\bdd t$ as the cone points and $m$ a vertical circle, and $t$ is rel
$\bdd t$ isotopic to a pair of arcs on $\bdd B$ of slope $p/q$.  See
[HT].  By definition $T(p/q)$ is equivalent to $T(p'/q')$ if $p/q
\equiv p'/q'$ mod $1$.  The tangle is trivial if $q = 0$ or $1$.  We
will always assume $q\geq 2$ and hence $T(p/q)$ is nontrivial.  Denote
by $E(p/q)$ the tangle space $E(t)$ if $(B,t,m) = T(p/q)$.

%2.2
\begin{defn} {\rm We use $\bar p = \bar p(p,q)$ to denote the mod $q$
    inverse of $-p$ with minimal absolute value, i.e., $\bar p$
    satisfies $p\bar p \equiv -1$ mod $q$, and $2|\bar p| \leq q$.
    Similarly, $\bar p_i$ denotes $\bar p(p_i, q_i)$ throughout the
    paper. }
\end{defn}

By a deformation of $B$ one can see that $T(p/q) = (B,t,m)$ can be
isotoped so that $t$ is rel $\bdd$ isotopic to a pair of vertical
arcs, and $m$ is a curve of slope $\bar p/q$ on $\bdd B$.  See Figure
2.1, where (a) is the standard picture of a $T(1/3)$, and (b) is
$T(1/3)$ after the deformation.  We have $\bar p = -1$ in this case,
so $m$ is a curve of slope $-1/3$ on $\bdd B$.  In general this can be
proved by lifting the boundary map to the universal cover of $\bdd B$,
considered as a sphere with four cone points of order 2.  It is
represented by a matrix $A$ with $\det (A) = 1$, which maps a line of
slope $p/q$ to $1/0$.  One can show that $A$ maps the line of slope
$1/0$ to $\bar p/q$.  We will use both points of view for $T(p/q)$ in
below.  

\bigskip
\leavevmode

\centerline{\epsfbox{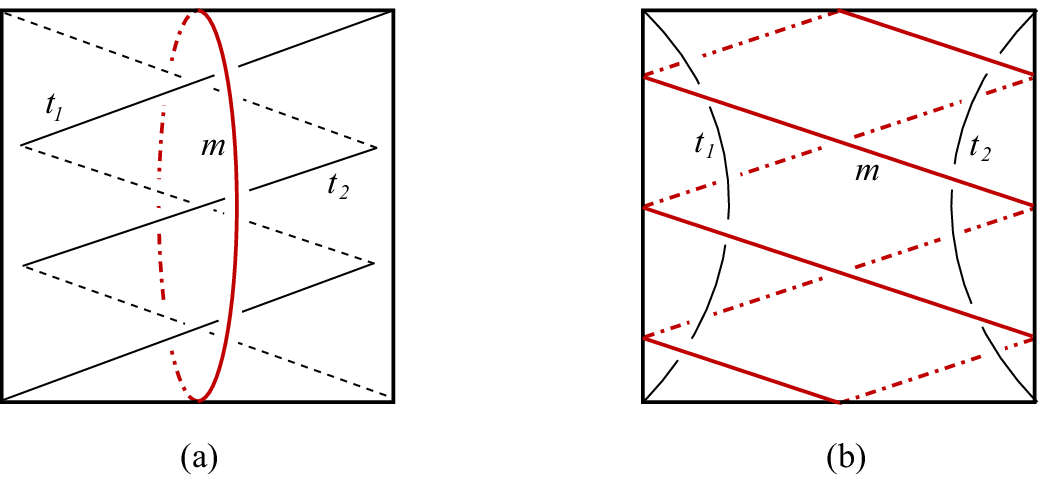}}
\bigskip
\centerline{Figure 2.1}
\bigskip

Let $E_0$ be an embedded disk in $B$ that separates the two arcs of
$t$ and intersects $m$ in $2q$ points, see Figure 2.2(a).  Let $E_i$
($i=1,2$) be embedded disks with $\bdd E_i = t_i \cup \alpha_i$, where
$\alpha_i$ is an arc on $\bdd B$ intersecting $m$ at $q$ points.  See
Figure 2.2(b).  These are chosen to be disjoint from each other.  Now
let $E_3$ be the disk in Figure 2.2(c), which intersects $E_0$ in one
arc but has interior disjoint from $E_1$ and $E_2$, and we have $\bdd
E_3 = t_1 \cup \b_1 \cup t_2 \cup \b_2$, where $\b_1, \b_2$ are arcs
on $\bdd B$, each intersecting $m$ at $|\bar p|$ points.  Thus $|\bdd
E_3 \cap m| = 2|\bar p|$.  Note that $E_3$ is unique up to isotopy
when $q>2$, and there are two such $E_3$ when $p/q = 1/2$.

\bigskip
\leavevmode

\centerline{\epsfbox{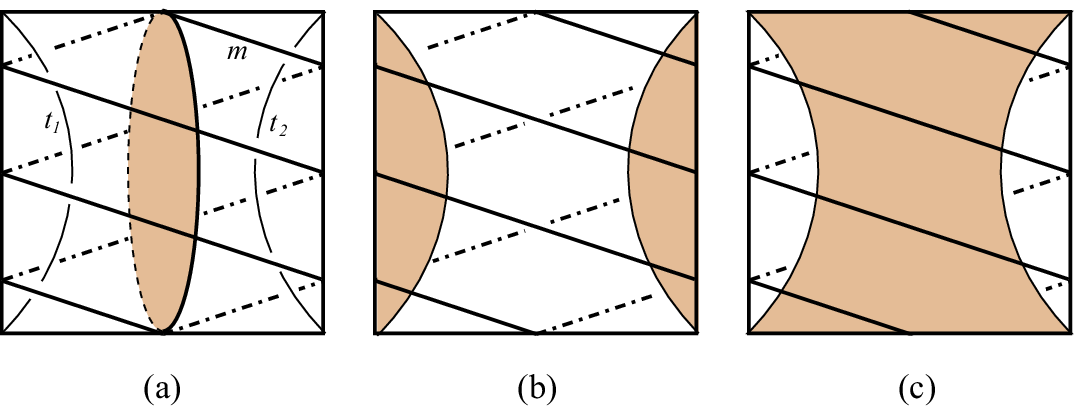}}
\bigskip
\centerline{Figure 2.2}
\bigskip

We call $E_0, E_1, E_2, E_3$ the {\it standard disks\/} for $(B, t)$.
We will also use the same notation $E_i$ to denote the disk $E_i \cap
E(t)$ when dealing with disks in $E(t)$.

An embedded disk $D$ in $E(t)$ is {\it tight\/} if it intersects $m$
and $A(t)$ minimally up to isotopy.  $D$ is an $(r,s)$ disk if $|\bdd
D \cap A(t)| = r$ and $|\bdd D \cap m| = s$.  See Definition 3.1 for
general definitions of tight curves and tight immersed disks.  We need
results which show that certain immersed disks are standard.  As a
warm up, we have the following lemma, which says that the standard
disks are the standard models of tight $(r,s)$ disks in $E(t)$ when
$r\leq 1$, or $r=2$ and $s$ is minimal.  The proof is a standard
innermost circle outermost arc argument by considering the
intersection of $D$ with $E_1 \cup E_2$, and is omitted.  Certain
version of (2) and (3) is also true for immersed disks when changing
isotopy to homotopy.  See Lemmas 3.5 and 3.6.

%2.3
\begin{lemma} Suppose $(B, t, m)$ is a $p/q$ rational tangle.  Let
  $\bar p$ be as in Definition 2.2.

  (1) If $D$ is am embedded tight $(0, s)$ disk in $E(t)$ then $s =
  2q$ and $D$ is $P$-isotopic to the standard disk $E_0$.

  (2) If $D$ is an embedded tight $(1, s)$ disk in $E(t)$ then $s = q$
  and $D$ is $P$-isotopic to the standard disk $E_1$ or $E_2$.

  (3) If $D$ is an embedded tight $(2, s)$ disk then $s \geq 2|\bar
  p|$, and if $s = 2|\bar p|$ then $D$ is $P$-isotopic to a standard
  disk $E_3$.

  (4) Standard disk of type $E_i$ are unique up to isotopy unless
  $q=2$ and $i=3$, in which case there are exactly two such isotopy
  classes. 
\end{lemma}

Denote by $|X|$ the number of components of $X$.  Note that if $X$ is
an immersed curve or surface in a manifold $M$ then it may have
intersection between different components, in which case $|X|$ denotes
the number of components before immersion, not the number of the
components of its image in $M$.  Thus two components of an immersed
surface which intersect in $M$ are still considered as different
components when counting $|X|$.

\section{Immersed disks in tangle spaces}

%3.1
\begin{defn} {\rm
  (1) A (possibly non-simple) closed curve $C$ on $\bdd E(t)$ is a
  {\em tight curve\/} if each component of $C\cap A(t)$, $C \cap
  P_1(t)$ and $C \cap P_2(t)$ is nontrivial.

  (2) A (possibly non-closed) curve $C$ is an {\em $(r,s)$-curve\/}
  (or a curve of type $(r,s)$) if $|\bdd D \cap A(t)| = r$, and $|\bdd
  D \cap m| = s$.  

  (3) A disk $D$ in $E(t)$ is a {\em $(r,s)$-disk\/} (or a disk of
  type $(r,s)$) if $\bdd D \cap \bdd E(t)$ is a $(r,s)$-curve on $\bdd
  E(t)$.  It is a {\em tight disk\/} if $\bdd D \subset \bdd E(t)$ is
  tight. }
\end{defn}

Note that if $C$ is a tight curve on $\bdd E(t)$ then each arc
component of $C\cap A(t)$ is homotopic to an embedded arc.  However,
$C\cap P_i(t)$ may contain arcs which have self intersections that
cannot be removed by homotopy.

%3.2
\begin{lemma} (1) Suppose $C$ is a tight curve on $F = \bdd E(t)$.
  Then it has minimal intersection with both $m$ and $A(t)$ up to
  homotopy.

  (2) Any curve $C$ on $F$ is homotopic to a tight curve $C'$, which
  is unique up to $P$-homotopy.  In particular, if two tight curves
  are homotopic then they are $P$-homotopic.

  (3) Suppose $C, C'$ are tight curves on $\bdd E(t)$.  If there are
  arcs $\alpha \subset C$ and $\beta \subset C'$ such that $\bdd
  \alpha = \bdd \beta$ and $\alpha \cup \beta$ is a trivial loop on
  $\bdd E(t)$, then $|\alpha \cap m| = |\beta \cap m|$ and $|\alpha
  \cap A(t)| = |\beta \cap A(t)|$.
\end{lemma}

\proof (1) If $C$ is homotopic to $C_1$ which has fewer intersection
with $m$, say, then the homotopy is a map $\varphi$ from an annulus $H
= S^1 \times I$ to $\bdd E(t)$.  Since $m \cup \bdd A(t)$ is embedded,
by transversality $\gamma = \varphi^{-1}(m \cup \bdd A(t))$ is an
embedded 1-manifold on $H$.  By a homotopy we may assume $\gamma$ has
no loops.  Since $|C \cap m| > |C_1 \cap m|$, $\gamma$ contains a
trivial arc on $H$ with both endpoints on $C$.  An outermost such arc
in $\gamma$ then cuts off a disk $D$ which gives rise to a homotopy
from an arc $\alpha$ of $C$ to an arc on $m$ or $\bdd A(t)$, and the
image of the interior of $D$ is disjoint from $m \cup \bdd A(t)$ and
hence is mapped into some $P_i(t)$ or $A(t)$.  It follows that
$\alpha$ is a trivial arc on $P_i(t)$ or $A(t)$, contradicting the
assumption that $C$ is tight.

(2) It is clear that any curve is homotopic to a tight curve.  We only
need to prove the uniqueness.  Let $\varphi: H \to \bdd E(t)$ be a
homotopy from $C_1$ to $C_2$ and assume $C_i$ are tight.  Consider
$\gamma = \varphi^{-1}(m \cup \bdd A(t))$.  As in the proof of (1) we
may assume $\varphi$ has no trivial loops, and the tightness of $C_i$
implies that there is no trivial arcs, hence each component of
$\gamma$ is an essential arc.  Deform $\varphi$ so that $\gamma$ is a
product $X \times I \subset S^1 \times I = H$, where $X$ is a finite
set in $S^1$.  Clearly $\varphi$ is now a $P$-homotopy since it maps
each ${z} \times I$ to an arc in some $P_j(t)$ or $A(t)$.

(3) Let $\varphi: D \to \bdd E(t)$ be a map with $\bdd D$ mapped to
$\alpha \cup \beta$, chosen to have minimal intersection with $m \cup
\bdd A(t)$.  Then $\varphi^{-1} (m \cup \bdd A(t))$ has no loops, and
each arc must have one end on $\alpha$ and the other on $\beta$ as
otherwise we can get a contradiction as above.  It follows that
$|\alpha \cap \gamma| = |\beta \cap \gamma|$ for each component
$\gamma$ of $m \cup \bdd A(t)$.  
\qed

%3.3
\begin{lemma}  
  Suppose that $(B,t, m)$ is a $p/q$ rational tangle with $q\geq 2$,
  and $D$ is a nontrivial $(0,s)$ disk in $E(t)$.  Then $s \geq 2q$.
\end{lemma}

\proof Let $E_0$ be the standard disk separating $t_1, t_2$, as
defined in Section 2.  Then $E_0$ cuts the tangle space $E(T)$ into
two solid tori $V_1, V_2$, and $\bdd E_0$ cuts the $4$-punctured disk
$\bdd B \cap E(t)$ into two twice punctured disks $Q_1, Q_2$.  We may
assume that $D$ is tight, and $D \cap E_0$ consists of arcs, each of
which is embedded in $E_0$.  

Assume to the contrary that $s<2q$.  Among all such disks, choose $D$
so that $k = |\bdd D \cap \bdd E_0|$ is minimal.  First assume $k>0$
and let $\alpha$ be an arc component of $D \cap E_0$, which cuts $E_0$
into $D'_1$ and $D'_2$.  Without loss of generality we may assume that
$|D'_1 \cap m| \leq q$.  Cutting $D$ along $\alpha$ and pasting two
copies of $D'_1$, we obtain two disks $D_1, D_2$.  Since $|\bdd D_1
\cap m| + |\bdd D_2 \cap m| = |\bdd D \cap m| + 2 |\bdd D'_1 \cap m| <
4q$, one of the $D_i$, say $D_1$, has $|\bdd D_1 \cap m| < 2q$.  Since
$D_1$ can be perturbed to have fewer intersection with $E_0$, by our
choice of $D$ the curve $\bdd D_1$ must be trivial.  Write $\bdd D_1 =
\a \cup \b$ with $\a \subset \bdd D$ and $\b \subset \bdd E_0$.  Since
both $D$ and $E_0$ are tight, by Lemma 3.2(3) we have $|\a \cap m| =
|\b \cap m|$.  We can now homotope $D$ so that $\a$ is deformed to
$\b$ and then push off $E_0$.  This reduces $|\bdd D \cap \bdd E_0|$
without changing $|\bdd D \cap m|$, contradicting the choice of $D$.

We now assume that $\bdd D \cap \bdd E_0 = \emptyset$, so $D \subset
V_1$, say.  Note that $m$ intersects $Q_1$ in $q$ arcs $\alpha_1, ...,
\alpha_q$, each cutting $Q_1$ into a surface $F_i \subset Q_1$ which
is the union of two longitudinal annuli.  If $\bdd D$ is disjoint from
some $\alpha_i$ then $\bdd D \subset F_i$, so $\bdd D$ would be null
homotopic on $F_i$ because $F_i$ is longitudinal while $\bdd D$ is
null homotopic in $V_1$, which contradicts the assumption that $D$ is
nontrivial.  Hence $|\bdd D \cap \alpha_i|>0$.  Let $\beta_i$ be an
arc on $E_0$ with $\bdd \beta_i = \bdd \alpha_i$.  If $|\bdd D \cap
\alpha_i|=1$ then we would have two closed curves $\bdd D$ and
$\alpha_i \cup \beta_i$ on the annulus $Q_1 \cup E_0$ intersecting
transversely at a single point, which is absurd.  Therefore $|\bdd D
\cap \alpha_i|\geq 2$ for each $i$, hence $s = |\bdd D \cap m| \geq
2q$.  \qed \bigskip

\bigskip
\leavevmode

\centerline{\epsfbox{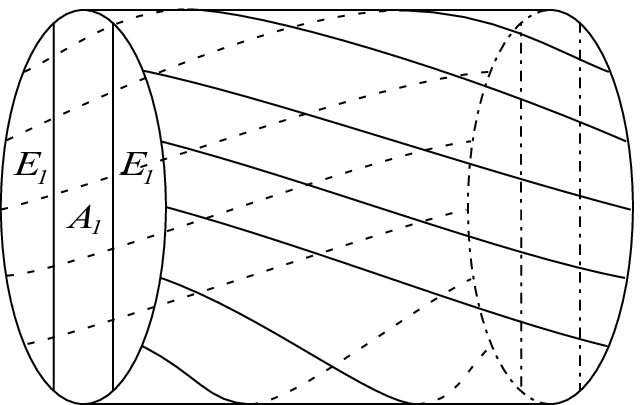}}
\bigskip
\centerline{Figure 3.1}
\bigskip

%3.4
\begin{lemma}
  Suppose that $(B,t,m)$ is a $p/q$ rational tangle.  Let $\bar p$ be
  as in Definition 2.2.  If $D$ is a nontrivial $(r,s)$ disk in
  $E(t)$, then $|s| \geq q$ if $r$ is odd, and $s \geq 2|\bar p|$ if
  $r$ is even.
\end{lemma}

\proof We may assume that $D$ is tight.  Consider the standard disks
$E_1, E_2$ defined in Section 2.  Let $E = E_1 \cup E_2$.  We proceed
by induction on (r, $|D\cap E|$).  By Lemma 3.3 the result is true if
$r = 0$.  Cutting $E(t)$ along $E$ produces a 3-ball on which $P(t)$
becomes an annulus $F$ and $m$ becomes a set of $2q$ essential arcs on
$F$.  See Figure 3.1 for the case $q=5$ and $\bar p = 2$.  Each
boundary component $\bdd_i$ of $F$ consists of 4 arcs, $\bdd_i = a'_i
\cup b'_i \cup a''_i \cup b''_i$, where $a'_i, a''_i$ are from $\bdd
A_i(t)$ and $b'_i, b''_i$ are copies of $\bdd E_i \cap P(t)$.  From
Figure 3.1 it is easy to see that the condition $2|\bar p| \leq q$ in
the definition of $\bar p$ implies that any arc $\gamma$ with
endpoints on $a'_1 \cup a''_1 \cup a'_2 \cup a''_2$ either is rel
$\bdd$ homotopic to an arc on one of the $a'_i$ or $a''_i$ and hence
is trivial on $P(t)$, or it intersects $m$ at least $|\bar p|$ times.
Moreover, if $\gamma$ has both endpoints on the same boundary
component of $F$ then it intersects $m$ at least $q$ times.  If $r$ is
odd then there is at least one arc with both endpoints on the same
component of $\bdd F$, which intersects $m$ at least $q$ times; if $r$
is even then there are at least two arcs, so they intersects $m$ at
least $2|\bar p|$ times.  Hence the result is true if $|D \cap E| =
0$.

Now suppose $|D \cap E| > 0$ and let $\alpha$ be an arc component of
$D\cap E_i$.  Then $\alpha$ cuts $D$ into $D', D''$ and it cuts $E_i$
into $E'_i, E''_i$, with $E''_i$ the one containing the arc on $A(t)$.
The four disks $D_1 = D' \cup E'_i$, $D_2 = D' \cup E''_i$, $D_3 = D''
\cup E'_i$, and $D_4 = D'' \cup E''_i$ are nontrivial.  For if $D_1$
is trivial, say, then we can homotope $D'$ to $E'_i$ and then push off
$E$ to reduce $|D \cap E|$; by Lemma 3.2(3) this will not change
$|D\cap m|$ and $|D \cap A(t)|$ and hence is a contradiction to the
minimality of $|D \cap E|$.

Suppose $D', D'', E'_i, E''_i$ are of types $(r_1, s_1)$, $(r_2,
s_2)$, $(0, s_3)$, $(1, s_4)$, respectively.  Then $r = r_1 + r_2$, $s
= s_1 + s_2$, and $s_3 + s_4 = q$.  The type of each $D_i$ is the sum
of the types of the corresponding subdisks, so the four disks $D_i$
are of the following types.

$D_1: \qquad (r_1,\,   s_1 + s_3)$

$D_2: \qquad (r_1+1,\, s_1 + s_4)$

$D_3: \qquad (r_2,\,   s_2 + s_3)$

$D_4: \qquad (r_2+1,\, s_2 + s_4)$

First assume $r$ is odd.  Then we may assume without loss of
generality that $r_1$ and $r_2+1$ are odd.  This implies that $r_1,
r_2+1 \leq r$.  Apply induction to $D_1$ and $D_4$ (note that they can
be deformed to have fewer intersection with $E$), we have $s_1 + s_3
\geq q$, $s_2 + s_4 \geq q$.  Adding these together gives $s + q = s_1
+ s_2 + s_3 + s_4 \geq 2q$, hence $s\geq q$.

Now assume $r$ is even.  If $r_1, r_2$ are odd then $r_i+1 \leq r$, so
by induction we have $s_1+s_3 \geq q$, $s_1 + s_4 \geq 2|\bar p|$,
$s_2 + s_3 \geq q$, and $s_2 + s_4 \geq 2|\bar p|$.  Adding these
together gives $2s + 2q = 2(s_1 + s_2 + s_3 + s_4) \geq 2q + 4|\bar
p|$, hence $s\geq 2|\bar p|$.  If both $r_i$ are even and nonzero then
a similar argument as above applies.  So we now assume that $r_1 = 0$
and $r_2 = r \geq 2$.  Then we can apply induction to $D_2$ and $D_3$
to get $s_1 + s_4 \geq q$ and $s_2 + s_3 \geq 2|\bar p|$.  Add these
together gives $q + s = s_1 + s_2 + s_3 + s_4 \geq q + 2|\bar p|$,
hence $s \geq 2|\bar p|$, as required.  \qed \bigskip

Suppose $F$ is am immersed surface represented by $\varphi: \tilde
F\to M$.  Let $S$ be an embedded surface in $M$.  Then $C =
\varphi^{-1}(S)$ is an embedded 1-manifold on $\tilde F$, and
$\varphi(C)$ is an immersed 1-manifold on $S$.  As in the embedded
case, we use $C = F \cap S$ to denote both the 1-manifold $C =
\varphi^{-1}(S)$ on $\tilde F$ and the immersed curve $\varphi: C \to
S$ on $S$.  To simplify notation we will not distinguish $\tilde F$
and $F$ and simply refer $C \subset \tilde F$ as $C \subset F$.  Thus
the curves $C$ above is still considered embedded on $F$ even though
it may have self intersection when considering $F$ as a subset of $M$.

The following two lemmas show that immersed tight $(1,q)$ disks and
$(2,2|\bar p|)$ disks are standard up to $P$-homotopy.

%3.5
\begin{lemma} Suppose $D$ is an immersed tight $(2,2|\bar p|)$ disk in
  $E(p/q)$.  Then it is $P$-homotopic to a standard disk $E_3$ defined
  in Section 2.
\end{lemma}

\proof By Lemma 3.2(2), if two tight curves are homotopic then they
are $P$-homotopic; hence we need only show that $D$ is homotopic to
some $E_3$.  

If $|D \cap E| = 0$ then $D$ lies on the 3-ball shown in Figure 3.1,
with $\bdd D$ consisting of one arc on each $A_i(t)$ and two arcs
$c_1, c_2$ on the annulus $F$ obtained by cutting $P(t)$ along $E$.
From Figure 3.1 we can see that the minimal intersection number
between $m$ and a nontrivial arc on $F$ is $|\bar p|$, and there are
exactly two such arc if $q>2$, or four such arcs if $q=2$.  These are
on the boundary of the standard disks of type $E_3$.  The assumption
that $D$ is a $(2,2|\bar p|)$ disk implies that $|c_i \cap m| = |\bar
p|$.  Hence we can deform $c_i$ by a homotopy of $(F,m)$ so that $\bdd
D$ matches the boundary of a standard disk $E_3$.  Since $E(t)$ is
irreducible, a further $P$-homotopy deforms $D$ to $E_3$.

Now assume $|D \cap E| > 0$.  Examine the proof of Lemma 3.4.  If one
of $\bdd D_i$ ($i=1,2,3,4$) is trivial then we may reduce $|D \cap E|$
by $P$-homotopy and the result follows by induction, so we may assume
that they are all nontrivial.

We may assume $s_1 \leq s_2$, so $s_1+s_2 = 2|\bar p| \leq q$ implies
$s_1 < q$.  Recall that $s_3 + s_4 = q$.  By the proof of Lemma 3.4,
if $r_1 = 0$ then $D_1$ would be a $(0, s_1 + s_3)$ disk, but since
$s_1 < q$ and $s_3 \leq q$, we would have $s_1 + s_3 < 2q$,
contradicting Lemma 3.3.  Therefore $r_1 = r_2 = 1$.  Since $D_1$ is a
$(1, s_1 + s_3)$ disk, by Lemma 3.4 we have $s_1 + s_3 \geq q$, so
$s_3 + s_4 = q$ implies $s_1 \geq s_4$.  Similarly we have $s_2 \geq
s_4$.  Since $D_2$ is of type $(2, s_1 + s_4)$, we have $2|\bar p| =
s_1 + s_2 \geq s_1 + s_4 \geq 2|\bar p|$, where the last inequality
follows from Lemma 3.4.  Hence we must have $s_2 = s_4$.  By symmetry
we have $s_1 = s_4$.  Thus both $D_2$ and $D_4$ are $(2, 2|\bar p|)$
disks, so by induction they are $P$-homotopic to $E_3$.  On the other
hand, by construction one of the $D_2, D_4$ has the property that it
has both arcs of $D_i \cap A(t)$ on the same component of $A(t)$,
which is a contradiction because $E_3$ has one edge on each of
$A_1(t)$ and $A_2(t)$.  \qed

%3.6
\begin{lemma} Suppose $D$ is a tight $(1,q)$ disk.  Then it is
  $P$-homotopic to the standard disk $E_1$ or $E_2$.
\end{lemma}

\proof Let $D_1, ..., D_4$ and $r_i, s_i$ be as defined in the proof
of Lemma 3.4.  The result is clear when $|D \cap E| = 0$, and it
follows by induction if one of the $D_i$ is trivial, so assume $|D
\cap E|>0$ and $\bdd D_i$ is nontrivial for $i=1,2,3,4$.  Since $r=
r_1 + r_2 = 1$, we may assume $r_1 = 0$.  The disk $D_1$ is of type
$(0,s_1+s_3)$, so by Lemma 3.3 we must have $s_1+s_3 \geq 2q$.  On the
other hand, we have $s_1+s_2 = s_3+s_4 = q$, so $s_2 = s_4 = 0$.  Now
$D_4$ is a $(2,0)$ disk, contradicting Lemma 3.4.  \qed
\bigskip

We note that a statement similar to Lemmas 3.5 and 3.6 is not true for
$(0,2q)$ disks.  A non-standard $(0,2q)$ disk can be formed by winding
$E_0$ around $A_i(t)$.  In other words, we can take $E_0$ and $n$
copies of $A_1$ or $A_2$ (but not both) and do cut and paste to form a
disk which is still an immersed $(0,2q)$ disk.  We call such a surface
an {\it $n$-winding disk\/} if $n$ is the minimal number of tubes
required in the above construction.  The curve $m$ cuts the boundary
of a $(0,2q)$ disk $D$ into $2q$ arcs.  We leave it to the reader to
verify that if $D$ is an $n$-winding disk then two of those $2q$ arcs
must have self intersection if $n>0$.

%3.7
\begin{lemma} Suppose $D$ is a tight $(0,2q)$ disk in $E(p/q)$.  Then
  $D$ is $P$-homotopic to an $n$-winding disk for some $n$.  In
  particular, if each component of $\bdd D \cap P(t)$ is an embedded
  arc then $D$ is $P$-homotopic to $E_0$.
\end{lemma}

\proof First assume $D \cap E = \emptyset$.  Cut $E(t)$ along $E$
produces a $D^2 \times I$, where $P(t)$ becomes an annulus $F$ and $m$
becomes a set of $2q$ parallel essential arcs on $F$, so it cuts $F$
into $2q$ squares.  See Figure 3.1.  Since $|\bdd D \cap m| = 2q$ and
$D$ is tight, $m$ cuts $\bdd D$ into $2q$ arcs, each of which is
nontrivial and hence connects one component of $m$ on $F$ to another.
It follows that each square contains exactly one arc of $\bdd D$,
which can be straightened inside of the square.  Therefore $\bdd D$ is
homotopic to the core of $F$, and $D$ is homotopic to $E_0$.  By Lemma
3.2(2) $D$ is $P$-homotopic to $E_0$.

Now assume $|D \cap E|>0$.  Consider an arc $\alpha \subset D \cap E$
which is outermost on $D$.  Let $D', D'', E', E''$ and $D_1, ..., D_4$
be as in the proof of Lemma 3.4, with $D'$ an outermost disk on $D$.
Then $D_1$ is a $(0,s_1 + s_3)$ disk and $D_3$ is a $(0, s_2+s_3)$
disk, so we have $s_1 + s_3 \geq 2q$ and $s_2 + s_3 \geq 2q$.  Since
$s_1 + s_2 = 2q$ and $s_3 \leq s_3+s_4 = q$, we must have $s_1 = s_2 =
s_3 = q$ and $s_4 = 0$.  It follows that there are exactly two
outermost disks on $D$, and all other components $Q$ of $D$ cut along
$D\cap E$ are bigons in the sense that it intersects $E$ in exactly
two arcs, and $Q \cap m = \emptyset$.  The union of $Q$ with two disks
on $E$ form a $(2,0)$ disk, so by Lemma 3.4 it must be trivial.  In
particular all components of $D\cap E$ are on the same disk $E_1$,
say.  It is now easy to see that each $Q$ is a tube (i.e.\ an annulus
parallel to $A_1(t)$) cut open, hence $D$ is an $n$-winding disk.  If
$n>0$ then by the above some component of $\bdd D \cap P(t)$ must have
self intersection.  Therefore if each component of $\bdd D \cap P(t)$
is embedded then we must have $n=0$ and hence $D$ is homotopic to
$E_0$.  \qed

\section{Immersed surfaces in essential position}

Recall that an embedded orientable surface $F$ in a 3-manifold is {\it
essential\/} if it is incompressible, $\bdd$-incompressible, and no
component of $F$ is boundary parallel.

%4.1
\begin{lemma}
  Let $F$ be an embedded orientable essential surface in a 3-manifold
  $M$.  Then no immersed essential arc $\alpha$ on $F$ is rel $\bdd
  \alpha$ homotopic to an arc on $\bdd M$.
\end{lemma}

\proof Suppose $\alpha$ is rel $\bdd \alpha$ homotopic to an arc
$\beta$ on $\bdd M$ and let $D$ be a null-homotopy disk bounded by $\a
\cup \b$.  If $M$ is reducible then since $F$ is essential we may
assume it is disjoint from a reducing sphere $S$, and we can then
modify $D$ if necessary to make it disjoint from $S$.  Therefore by
decomposing along $S$ we may assume that $M$ is irreducible.
Similarly, since $F$ is incompressible and $\bdd$-incompressible and
$M$ is irreducible, $F$ can be isotoped and $D$ can be modified to be
disjoint from any $\bdd$-reducing sphere, hence by cutting along
$\bdd$-reducing disks if necessary we may also assume that $M$ is
$\bdd$-irreducible.

Let $2M$ be the double of $M$ along $\bdd M$, and let $2F$ be the
double of $F$ along $\bdd F = F\cap \bdd M$.  Then the double of
$\alpha$ is an essential curve on $2F$ which is null homotopic in
$2M$.  Thus $2F$ is not $\pi_1$-injective and hence is compressible.
Let $D$ is a compressing disk of $2F$ in $2M$ such that $|D \cap \bdd
M|$ is minimal.  Since $F$ and $\bdd M$ are incompressible in $M$, by
an innermost circle argument one can show that $D\cap \bdd M$ has no
circle component, and we must have $D\cap \bdd M \neq \emptyset$
because $D$ cannot be a compressing disk of $F$.  An outermost arc of
$D\cap \bdd M$ then cuts off a $\bdd$-compressing disk of $F$ in $M$.
\qed

%4.2
\begin{lemma}
  Let $F$ be an embedded orientable essential surface in a 3-manifold
  $M$.  Suppose $\rho: \tilde M \to M$ is a finite cover.  Then the
  surface $\tilde F = \rho^{-1}(F)$ is essential in $\tilde M$.
\end{lemma}

\proof Let $\tilde F_1$ be a component of $\tilde F$ and let $F_1$ be
the corresponding component in $F$.  If $\tilde F_1$ is boundary
parallel then it cuts off a regular neighborhood of a boundary
component $\tilde T$ of $\tilde M$, which projects to a regular
neighborhood of a boundary component $T$ of $M$, hence $F_1$ would
also be boundary parallel, contradicting the assumption that $F$ is
essential.  

A compressing disk of $\tilde F_1$ would map to an immersed disk in
$M$ with boundary an essential curve on $F_1$, which contradicts the
fact that an incompressible surface is $\pi_1$-injective.  Therefore
$\tilde F$ is incompressible in $\tilde M$.  

If $\tilde F$ is $\bdd$-compressible then a $\bdd$-compressing disk of
$\tilde F$ in $\tilde M$ projects to a disk in $M$ whose boundary
consists of an essential arc on $F$ and an arc on $\bdd M$, which
contradicts Lemma 4.1.
\qed

%4.3
\begin{defn} {\rm Suppose $L$ is a link in a 3-manifold $M$ and $F$ is an
  embedded essential surface in $E(L) = M - \Int N(L)$ with nonempty
  non-meridional boundary slope on each boundary component of $E(L)$.
  Let $\hat S$ be an immersed surface in $M$, and let $S = \hat S \cap
  E(L)$.  Then $\hat S$ is said to be {\em in essential position with
  respect to $F$\/} if $\bdd S \cap L \neq \emptyset$, each
  arc component of $S \cap F$ is essential on both $S$ and $F$, and
  each circle component is nontrivial on $S$.}
\end{defn}

The following lemma is essentially due to Gabai [Ga].  

%4.4
\begin{lemma}[Thin Position Lemma] Suppose $L$ is a link in
  $S^3$, and $F$ an embedded essential surface in $E(L) = S^3 - \Int
  N(L)$ with nonempty non-meridional boundary slopes on each boundary
  component of $E(L)$.  Then there exist an embedded sphere $\hat S$
  in $S^3$ which is in essential position with respect to $F$.
\end{lemma}

\proof When $L$ is a knot this is [Ga, Lemma 4.4].  If $L$ splits, one
can use the incompressibility of $F$ to find a splitting sphere
disjoint from $F$, and proceed by induction to find $\hat S$ on the
side containing $F$; hence we may assume that $L$ is non-split.  In
this case the proof is the same as in [Ga, Lemma 4.4].  Put $L$ in
thin position and let $\hat S$ be a thin level surface.  Since $L$ is
non-split, $\hat S$ has nonempty intersection with $L$.  Isotoping $F$
properly and using the argument in [Ga, Lemma 4.4] one can show that
$\hat S$ can be chosen so that it has neither high disk nor low disk,
in which case the arc components of $S \cap F$ are essential on both
$F$ and $S$.  Since $F$ is incompressible, any circle of $F\cap S$
which is trivial on $S$ bounds a disk on $F$, so one can modify $S$ by
cut and paste to get rid of those curves.  We refer the reader to the
proof of [Ga, Lemma 4.4] for details.  \qed

The following can be considered as an immersed version of the thin
position lemma.  It works for links in manifolds with finite
fundamental group.

%4.5
\begin{lemma}[Immersed Thin Position Lemma]
  Suppose $M$ is a closed irreducible 3-manifold with finite
  fundamental group, $L$ a link in $M$, and $F$ an embedded essential
  surface in $E(L) = M - \Int N(L)$ with nonempty non-meridional
  boundary slopes on each boundary component of $E(L)$.  Then there
  exists an immersed sphere $Q'$ in $M$ which is in essential position
  with respect to $F$.
\end{lemma}

\proof The universal cover $\tilde M$ of $M$ is a simply connected
closed 3-manifold and hence by the Poincare conjecture proved by
Perelman [Pr], it is an $S^3$.  The lifting of $L$ is a link $\tilde
L$ in $\tilde M$, and the lifting of $F$ is a surface $\tilde F$ in
$E(\tilde L) = S^3 - \Int N(\tilde L)$.  Since $F$ is an essential
embedded surface in $E(L)$, by Lemma 4.2 $\tilde F$ is an essential
embedded surface in $\tilde M = E(\tilde L)$.  By the Thin Position
Lemma above, there is a sphere $\tilde Q'$ in $S^3$ which is in
essential position with respect to $\tilde F$.  Let $Q'$ be the
projections of $\tilde Q'$ in $M$, let $\tilde Q = \tilde Q' \cap
E(\tilde L)$, and let $Q = Q' \cap E(L)$.  The preimage of an arc in
$Q\cap F$ may consist of several arcs, but at least one of them is in
$\tilde Q \cap \tilde F$, and an inessential arc would lift to an
inessential arc.  Therefore all arcs of $Q\cap F$ are essential on
both $F$ and $Q$.  Similarly all circle components of $Q \cap F$ are
essential on $Q$.  \qed

\section{Elementary surfaces}

We will always use $K_r$ to denote the core of the Dehn filling solid
torus in the surgered manifold $K(r)$.  Given an immersed surface $F$
in $E(K\cup C)$ with each boundary component either a meridional curve
on $\bdd N(C)$ or a curve of slope $r$ on $\bdd N(K)$, denote by $\hat
F$ the closed surface obtained by attaching a meridian disk of $N(K_r
\cup C)$ to each boundary component of $F$.

Suppose $K = K(r_1, r_2, r_3)$ is a Montesinos knot of length 3, where
$r_i = p_i/q_i$ and $q_i \geq 2$.  Let $L = K \cup C$, where $C$ is
the axis of $K$ as shown in Figure 5.1.  Let $E(L) = S^3 - \Int N(L)$
be the exterior of $L$.  Denote by $V$ the solid torus $S^3 - \Int
N(C)$.  There are three disks $D_1, D_2, D_3$ cutting the pair $(V,
K)$ into three rational tangles $(B_i, t_i, m_i)$ of slope $r_i$.  Let
$P_i$ be the twice punctured disk $D_i \cap E(L)$, and let $P = P_1
\cup P_2 \cup P_3$.  Then $P = P_1 \cup P_2 \cup P_3$ cuts $E(K\cup
C)$ into three tangle spaces $E(t_1)$, $E(t_2)$, $E(t_3)$.  The
surface $\bdd E(t_i)$ is the union of $P_i, P_{i+1}$ and three annuli
$A_0(t_i)$, $A_1(t_i)$, $A_2(t_i)$, where $A_0(t_i) = \bdd E(t_i) \cap
\bdd N(C)$ and $A_1(t_i) \cup A_2(t_i)$ are the ones on $\bdd N(K)$.

\bigskip
\leavevmode

\centerline{\epsfbox{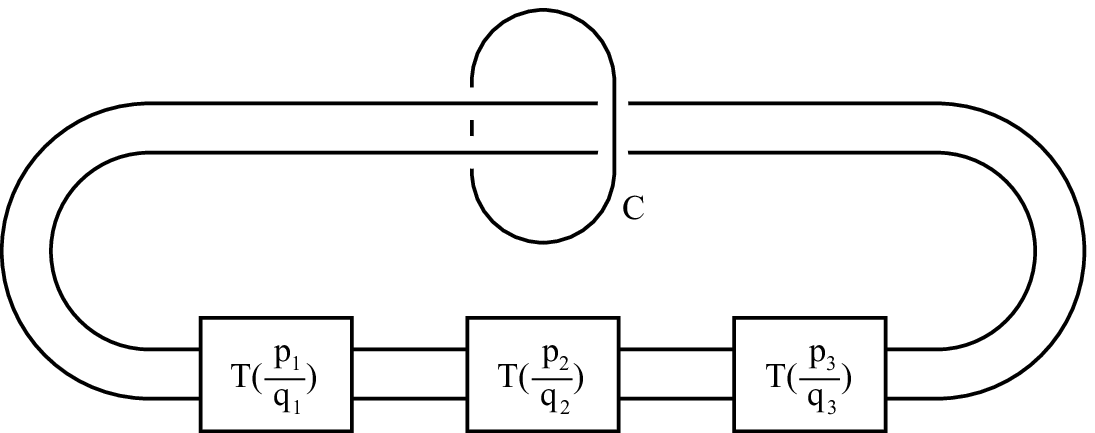}}
\bigskip
\centerline{Figure 5.1}
\bigskip

%5.1
\begin{defn} {\rm 
  An immersed surface $F$ in $E(K\cup C)$ is an {\em elementary
    surface\/} if the following holds.

  (1) Each component of $F \cap E(t_i)$ is a tight disk of type
  $(0,2q_i)$, $(1, q_i)$ or $(2, 2|\bar p_i|)$;

  (2) at least one $F \cap E(t_i)$ has no disk of type $(0, 2q_i)$;

  (3) if $q_i = 2$ then all disks of type $(2, 2|\bar p_i|)$ in
  $E(t_i)$ are homotopic to each other; and

  (4) type $(0, 2q_i)$ and type $(2, 2|\bar p_i|)$ disks do not appear
  simultaneously in any $E(t_i)$. }
\end{defn}

%5.2
\begin{lemma} Suppose $K$ is hyperbolic, and there is an {\em
    embedded} closed torus or Klein bottle $\hat F$ in $K(r)$ such
  that $F = \hat F \cap E(K \cup C)$ is an elementary surface.  Then
  $K(r)$ is toroidal.
\end{lemma}

\proof By considering the boundary of a regular neighborhood of $F$ if
necessary, we may assume $F$ is orientable and hence is a punctured
torus.   Since $F$ is embedded
and $F \cap E(t_i)$ is the union of standard disks in $E(t_i)$ for all
$i$, the surface $\hat F$ is a candidate surface as defined in [HO].
Since $\hat F$ is a torus, by [Wu3, Lemma 7.1] the knot $K$ and the
boundary slope $r$ of $\hat F$ are among those listed in [Wu3, Theorem
1.1].  By [Wu3, Theorem 1.2] $\hat F$ is incompressible and hence $K(r)$
is toroidal.  \qed
\bigskip

The main result of this section is Proposition 5.6, which shows that
the above lemme is still true for immersed $\hat F$.  

%5.3
\begin{lemma} Suppose $F$ is an immersed elementary surface in $E(K \cup C)$.

  (1) For each $i$, $F$ can be deformed so that $F\cap E(t_i)$ is
  embedded.  

  (2) $F$ can be deformed so that $F \cap P$ is embedded.
\end{lemma}

\proof By Lemmas 3.5--3.7 each component of $F \cap E(t_i)$ can be
deformed to a standard disk or an $n$-winding disk.  By condition (2)
in Definition 5.1, up to relabeling we may assume that there is no
disk of type $(0,2q_1)$ and hence no $n$-winding disk in $E(t_1)$ for
$n>0$.  Since standard disks are embedded, each arc of $F \cap P_1$
and $F \cap P_2$ is embedded up to homotopy.  Since $K$ is a knot, at
most one of the $q_j$ is even, hence without loss of generality we may
assume that $q_2$ is odd.  It means that a strand of $t_2$ has one
endpoint on each of $P_2$ and $P_3$.  If $F \cap E(t_2)$ contains an
$n$-winding disk $D$ for some $n>0$ then $D \cap P_2$ would be an arc
which is not homotopic to an embedded arc, a contradiction.  Therefore
$E(t_2)$ has no such disk, which implies that each component of $F
\cap E(t_2)$ is also homotopic to a standard disk.  Similarly, since
each arc in $F\cap (P_1 \cup P_3)$ is homotopic to an embedded arc,
there is no $n$-winding disk in $E(t_3)$ for any $n>0$.  Therefore
each component of $F\cap E(t_j)$ is homotopic to a standard disk for
all $j$.

By definition, standard disks can be homotoped in $E(t_j)$ to be
disjoint from each other except that (i) a disk of type $E_0$
intersects a disk of type $E_3$ essentially, and (ii) when $q_j=2$
there are two possible homotopy classes of type $E_2$ which intersect
each other essentially.  However, these have been excluded by
conditions (3) and (4) in Definition 5.1.  Therefore $F\cap E(t_j)$
can be deformed to be a set of mutually disjoint standard disks.  This
proves (1).

It may not be possible to do the above simultaneously for all the
three tangles, but it implies that each $F\cap P_i$ can be deformed to
be embedded, hence we can homotope $F$ so that $F \cap P$ is embedded,
as stated in (2).  \qed

%5.4
\begin{lemma} Let $P$ be a twice punctured disk, and $P_0, P_1$
  the two boundary copies of $P$ in $P \times I$.  Suppose $D_1$ is an
  embedded disk in $P\times I$ intersecting each $P_i$ at a single
  essential arc $a_i$.  Let $D_2$ be an immersed disk in $P\times I$
  with $b_i = D_2 \cap P_i = a_i$ for $i=0,1$.  Then $D_2$ is
  homotopic to $D_1$ rel $b_0 \cup b_1$.
\end{lemma}

\proof By an isotopy of $P \times I$ we may assume that $D_1$ is a
product disk $\alpha \times I$.  Let $\b_1, \b_2$ be two arcs on $P$
such that $\a, \b_1, \b_2$ are mutually disjoint, mutually
non-parallel, and $\b_1 \cup \b_2$ cuts $P$ into a disk.  We may write
$\bdd D_1 = a_0 \cup a_2 \cup a_1 \cup a_3$ and $\bdd D_2 = b_0 \cup
b_2 \cup b_1 \cup b_3$, so that $\bdd a_2 = \bdd b_2$, $\bdd a_3 =
\bdd b_3$, and $a_2, a_3, b_2, b_3$ are arcs on $\bdd P \times I$.
Put $Q_i = \b_i \times I \subset P \times I$ and let $Q = Q_1 \cup
Q_2$.  Homotope $D_2$ rel $b_0 \cup b_1$ so that it has minimal
intersection with $Q$.  Let $\varphi: D \to D_2$ be the immersion.
Then $\varphi^{-1}(Q)$ is an embedded 1-manifold on $D$.  Using the
fact that the arc $\a, \b_1, \b_2$ are essential mutually
non-homotopic proper arcs in $P$ one can apply an innermost circle
outermost arc argument to show that $D_2 \cap Q = \emptyset$.  Since
$Q$ cuts $\bdd P \times I$ into disks, we see that $a_2 \cup b_2$ and
$a_3 \cup b_3$ are trivial loops, hence up to homotopy rel $b_0 \cup
b_1$ we can deform $\bdd D_2$ to $\bdd D_1$.  Since $P \times I$ is
irreducible, one can further deform $D_2$ to $D_1$ by a homotopy rel
$\bdd D_2$.  \qed \bigskip

We now consider a set of {\it oriented\/} loops $C$ on a torus $T$
with oriented meridian-longitude pair $(\mu, \lambda)$.  Then $C$
represents some $a\mu + b \lambda$ in $H_1(T)$.  Define $b/a$ to be
the {\em slope\/} of $C$.  Note that $C$ is not required to be
embedded.  If $C$ has a crossing as in Figure 5.2(a), we can change it
to that in Figure 5.2(b).  This operation is called {\em smoothing\/}
a crossing.

\bigskip
\leavevmode

\centerline{\epsfbox{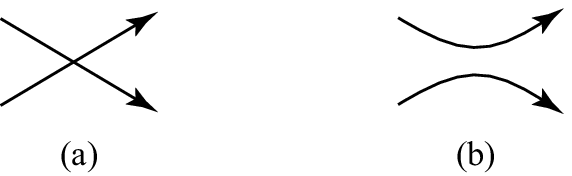}}
\bigskip
\centerline{Figure 5.2}
\bigskip

%5.5
\begin{lemma}
  Let $C$ be a set of oriented immersed curves on a torus $T$.  Let
  $C'$ be obtained from $C$ by smoothing crossings.  Then $C$ and $C'$
  have the same slope on $T$.
\end{lemma}

\proof  This is obvious since the curves in Figure 5.2(a) and 5.2(b)
represents the same homology class in $H_1(T)$.
\qed

%5.6
\begin{prop} Suppose $K = K(p_1/q_1,\, p_2/q_2,\, p_3/q_3)$ is
  hyperbolic and $\hat F$ is an immersed $\pi_1$-injective torus in
  $K(r)$.  If $F = \hat F \cap E(K\cup C)$ is an elementary surface,
  then $K(r)$ is toroidal.
\end{prop}

\proof As before, let $P = P_1 \cup P_2 \cup P_3$ be the three twice
punctured disks that cut $E(K\cup C)$ into $E(t_1) \cup E(t_2) \cup
E(t_3)$.  Let $W_i = P_i \times I$ be a collar of $P_i$.  Then we can
rewrite $E(K\cup C)$ as the union $(\cup W_i) \cup_R (\cup E(t_i))$,
where $W_i$ and $E(t_j)$ have mutually disjoint interiors, and $R =
(\cup W_i) \cap (\cup E(t_i))$ is the union of 6 twice punctured disks.

Since $F$ is elementary, by Lemma 5.3(2) we may assume that $F \cap P$
is embedded, so $F \cap W_i$ consists of mutually disjoint product
disks.  Since the $E(t_i)$ are mutually disjoint, by Lemma 5.3(1), we
can deform $F_i = F \cap E(t_i)$ in $E(t_i)$ so that it is embedded.
Note that this can be done by a homotopy of the pair $(E(t_i), R \cap
E(t_i))$, so they can be combined and then extended into a small
neighborhood of $R$ in $\cup W_j$ to form a homotopy of $F$.  After
this homotopy $F \cap W_i$ may no longer be embedded; however, since
we started with product disks and the homotopy maps each component of
$R$ to itself, each component of $F \cap W_i$ is still homotopic to a
product disk.

\bigskip
\leavevmode

\centerline{\epsfbox{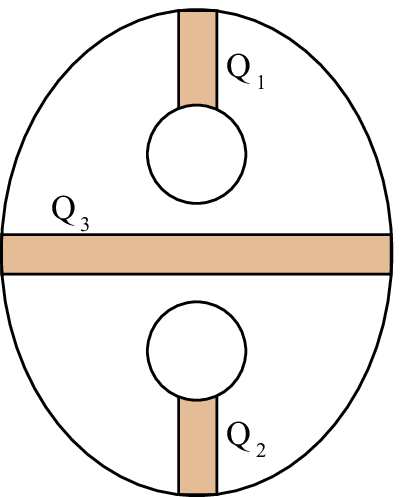}}
\bigskip
\centerline{Figure 5.3}
\bigskip

Let $P'_i, P''_i$ be the two copies of $P_i$ on $\bdd W_i$.  Since $F
\cap E(t_i)$ is a union of mutually disjoint standard disks, we may
assume that $F \cap P'_i$ consists of embedded essential arcs in the
small disks $Q_1 \cup Q_2\cup Q_3$ shown in Figure 5.3.  Since each
component $D$ of $F\cap W_i$ is a disk, the two arcs $D\cap P'_i$ and
$D\cap P''_i$ must lie in $Q_j \times 0$ and $Q_j\times 1$ for the
same $j$ because two essential arcs in different $Q_j$ are not
homotopic to each other.  Let $E$ be an embedded disk in $Q_j \times
I$ such that $E \cap P'_i = D \cap P'_i$ and $E \cap P''_i = D \cap
P''_i$.  By Lemma 5.4 $D$ is rel $D \cap (P'_i \cup P''_i)$ homotopic
to $E$.  Therefore up to homotopy we may assume that each component
$D$ of $F \cap W_i$ is a flat embedded disk in $Q_j \times I$, so
different components $D', D''$ of $F \cap W_i$ are either disjoint or
intersect at a single arc, as shown in Figure 5.4(a).

\bigskip
\leavevmode

\centerline{\epsfbox{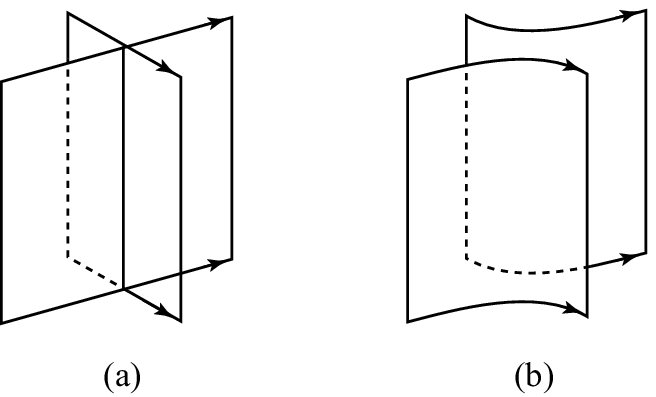}}
\bigskip
\centerline{Figure 5.4}
\bigskip

Orient $\bdd F$ so that all arcs of $\bdd F \cap W_i$ runs
monotonically from $P'_i$ to $P''_i$.  We can then perform a double
edge smoothing of $D', D''$ to obtain two new disks $D'_1, D''_1$, as
shown in Figure 5.4(b).  After doing this for all double curves, the
surface $F \cap W_i$ becomes a set of mutually disjoint embedded
disks, and the original surface $F$ now becomes an embedded surface
$F'$ in $E(K\cup C)$.  By definition $F'$ is an elementary
surface.  Note that the double curve smoothing does not change the
Euler characteristic of a surface.  By Lemma 5.5 $F'$ has the same
boundary slope as that of $F$.  Since both $\bdd F'$ and $\bdd F$ are
embedded, we see that they have the same number of boundary curves.
Let $\hat F'$ be the closed surface obtained by adding meridian disks
of $N(K_r\cup C)$ to $\bdd F'$.  The above implies that $\chi(\hat F') =
\chi(\hat F) = 0$, hence $\hat F'$ is a torus or Klein bottle, and it
is embedded.  It now follows from Lemma 5.2 that $K(r)$ is toroidal.
\qed

\section{Intersection graphs}

We continue using the notations in Section 5, so $K = K(p_1/q_1,\,
p_2/q_2,\, p_3/q_3)$, $L = K \cup C$, and $P = P_1 \cup P_2 \cup P_3$,
cutting $E(L)$ into three tangle spaces $E(t_1), E(t_2), E(t_3)$.  Let
$K(r)$ be the manifold obtained by Dehn surgery along a slope $r$ on
$K$.  Denote by $K_r$ the core of the Dehn filling solid torus in
$K(r)$, and let $L_r = K_r \cup C$.

%6.1
\begin{lemma} $P$ is an essential surface in $E(L)$.
\end{lemma}

\proof If $P$ is compressible then the compressing disk would separate
the two strings of $t$ and $t$ would be a trivial tangle,
contradicting the assumption that $q_i \geq 2$.  It is well known that
if $\bdd M$ is a set of tori then any connected incompressible surface
in $M$ is also $\bdd$-incompressible unless it is an annulus.  Hence
$P$ is also $\bdd$-incompressible.  Since each component of $P$ has
nonempty boundary, this also implies that no component of $P$ is
boundary parallel.  \qed

%6.2
\begin{lemma} Suppose $M = K(r)$ is Seifert fibered.  Then there is an
  immersed torus or sphere $\hat Q$ in $M$ which is in essential
  position with respect to $P$.
\end{lemma}

\proof By Lemma 6.1 $P$ is an embedded essential surface in $E(L_r) =
E(L)$.  Therefore the result follows from Lemma 4.5 if
$\pi_1(K(r))$ is finite, and the surface $\hat Q$ is a sphere in this
case.

Now assume $\pi_1(K(r))$ is infinite.  By [Sc] a Seifert fibered
manifold admits one of six 3-manifold geometries of Thurston.  Since
$K(r)$ is irreducible by [Wu1], it does not admit $S^2 \times R$
geometry, and since $\pi_1(K(r))$ is infinite it does not admit $S^3$
geometry.  Hence the orbifold $X$ of $K(r)$ must be euclidean or
hyperbolic and therefore has infinite orbifold fundamental group as
defined in [Sc].  The torus $\hat Q$ in $K(r)$ that projects to a
curve on $X$ with infinite order in the orbifold fundamental group is
then an immersed torus in $K(r)$ which is $\pi_1$-injective in $K(r)$.

Let $Q = \hat Q \cap E(L_r)$.  Choose $\hat Q$ to have minimal
intersection with $L_r$ and then homotoped so that $(|\bdd Q \cap \bdd
P|, |Q \cap P|)$ is minimal in lexicographic order.  In particular,
$\bdd Q$ intersects $\bdd P$ minimally.  We need to show that $\hat Q$
is in essential position with respect to $P$.  Since $P$ is essential,
we may get rid of loops in $Q\cap P$ which is trivial on $Q$.  We need
to show that all arc components $\alpha$ on $P\cap Q$ are nontrivial
on both $P$ and $Q$.  Note that $\alpha$ is embedded on $Q$ but may be
immersed on $P$.

Assume $\alpha$ is inessential on both $P$ and $Q$.  Then it is rel
$\bdd \alpha$ homotopic to an arc $\beta$ on $\bdd Q$ and an arc
$\gamma$ on $\bdd P$, so the incompressibility of $\bdd N(K)$ implies
that $\b$ is homotopic to $\gamma$ and hence there is a homotopy to
reduce $|\bdd Q \cap \bdd P|$, contradicting its minimality.

Now assume $\alpha$ is inessential on $P$ but essential on $Q$.  Since
$\alpha$ is trivial on $P$, we may deform $Q$ near $P$ to make
$\alpha$ embedded on $P$.  One can then push $\hat Q$ through the disk
on $P$ cut off by $\alpha$ to reduce $|\hat Q \cap L_r|$ by 2,
contradicting its minimality.

If $\alpha$ is essential on $P$ but inessential on $Q$ then it cuts
off a disk $D$ on $Q$ and hence is rel $\bdd \alpha$ homotopic to
the arc $\bdd D - \Int \alpha \subset \bdd E(L_r)$.  Since $P$ is
essential by Lemma 6.1, this contradicts Lemma 4.1.
\qed
\bigskip

Let $\hat Q$ be as in Lemma 6.2 and put $Q = \hat Q \cap E(L_r)$.
Since $P$ is essential and embedded, $P \cap Q$ is an embedded compact
1-manifold in $Q$.  However it may not be embedded in $P$ because $Q$
is immersed in $E(L)$.  Since $\hat Q$ is in essential position with
respect to $P$, no component of $Q\cap P$ is a trivial circle or
trivial arc on $Q$.

Define a {\it generalized graph\/} to be a graph with possibly some
valence 2 vertices removed from the vertex list.  Thus it is a graph
except that some loop components may not have vertices on it.

%6.3
\begin{defn}  {\rm
  Consider the disks of $\hat Q - \Int Q = \hat Q \cap N(L_r)$ as fat
  vertices, the arc components of $P\cap Q$ as edges, and the
  circles of $P\cap Q$ as loops.  This produces a generalized graph
  $G$ on $\hat Q$, called the {\em intersection graph\/} of $\hat Q$
  and $P$. }
\end{defn}

There are two types of vertices on $G$.  A component of $\hat Q \cap
N(C)$ is called a {\it small vertex}, and a component of $\hat Q \cap
N(K)$ is called a {\it large vertex}.  Since a meridian of $C$
intersects $P$ at three points, one on each $P_i$, we see that each
small vertex has valence 3.  Denote by $\Delta = \Delta(\mu, r)$ the
minimal intersection number between a meridian $\mu$ of $K$ and the
surgery slope $r$.  Then the surgery slope $r$ intersects $\bdd P$
minimally at $6\Delta$ points, so each large vertex has valence
$6\Delta$, with $2\Delta$ endpoints on each $P_i$.  In particular, if
$r$ is an integer slope, we have $\Delta =1$, in which case the
valence of each large vertex is $6$.

Write the boundary of the tangle space $E(t_i)$ as $\bdd E(t_i) = P_i
\cup P_{i+1} \cup A_i \cup A'_i \cup A''_i$, where $P_i$ are the twice
punctured disk $D_i \cap E(L)$ defined above, $A_i = \bdd N(C) \cap
E(t_i)$, and $A'_i, A''_i$ are the two annuli $\bdd N(K) \cap E(t_i)$.
Since $\hat Q$ is in essential position with $P$, we have the
following lemma.

%6.4
\begin{lemma} 
  Let $\sigma$ be a face of $G$ lying in $E(t_i)$.  Then $\bdd \sigma$
  intersects each of the above five subsurfaces of $\bdd E(t_i)$ in
  essential arcs and hence is a set of tight curves.  In particular,
  each disk face of $G$ is an essential disk in some $E(t_i)$.
\end{lemma}

An arc of a face $\sigma$ on the boundary of a fat vertex $v$ is
called a {\it corners\/} of $\sigma$ at $v$.  Note that when shrinking
each fat vertex to a single point a corner becomes a vertex on the
boundary of the face $\sigma$.  A corner at $v$ is {\it large\/} or
{\it small\/} according to whether $v$ is a large vertex or a small
vertex.  Thus a large corner lies on $A'_i$ or $A''_i$ for some $i$
while a small corner is on the other annulus $A_i$ and hence
intersects $m_i$ at a single point.  Therefore a disk face $\sigma$ of
$G$ is of type $(r,s)$ if and only if it has $r$ large corners and $s$
small corners.  The results in Section 3 now apply to the disk faces
of $G$.  In particular, an $(r, s)$ face in $E(t_i)$ has $s \geq 2q_i$
if $r=0$, $s\geq q_i$ if $r$ is odd, and $s\geq 2 |\bar p_i|$ if $r$
is even.

\section{Euler number of an angled surface}

An {\it angled surface\/} is a compact surface $\sigma$ with a set of
points $V = (v_1, ..., v_n)$ on $\bdd \sigma$ called vertices or
corners, and an angle $\a_i$ assigned to each corner $v_i$, with $0
\leq \a_i < \pi$.  When $\sigma$ is a disk, it is a polygon or
$n$-gon.  We will always use $\bar \a$ to denote the external angle
$\bar \a = \pi - \a$.

%7.1
\begin{defn} {\rm The {\em angled Euler number\/} of an angled
    surface $\sigma$ with corner angles $\a_1, ..., \a_n$ is defined as
    $$e(\sigma) = \chi(\sigma) - \frac 1{2\pi} \sum \bar \a_i$$ }
\end{defn}

Recall that a generalized graph is a graph on which some of the loops
may not have vertices.

%7.2
\begin{lemma} Let $G$ be a generalized graph on a closed surface $F$,
  cutting it into angled surfaces $\sigma_1, ..., \sigma_m$, such that
  the sum of angles around each vertex is at least $2\pi$.  Then $
  \chi(F) \leq \sum e(\sigma_i)$, and equality holds if and only if
  the sum of angles around each vertex is $2\pi$.
\end{lemma}

\proof Note that adding vertices to loops with an angle $\pi$
at each corner of the new vertices will not change the angled Euler
number of the faces.  Therefore by adding such vertices if necessary
we may assume that $G$ is a genuine graph.  We assume that the sum of
the angles around each vertex of $G$ is exactly $2\pi$.  The proof for
the other case is similar.

Let $n_i$ be the number of vertices on $\bdd \sigma_i$, let $v_{ij}$ be the
vertices on $\bdd \sigma_i$, and let $\alpha_{ij}$ be the angles at
$v_{ij}$. Note that $n_i$ is also the number of edges on $\bdd \sigma_i$.
Denote by $E$ and $V$ the number of edges and vertices of $G$,
respectively.  Then
\begin{eqnarray*}
\sum_i e(\sigma_i) &=& \sum_i [\chi(\sigma_i) - \frac 1{2\pi} \sum_j (\pi -
\a_{ij})] \\
&=& \sum_i [\chi(\sigma_i) - \frac 12 n_i + \sum_j \frac{\a_{ij}}{\pi}\;] \\
&=& \sum \chi(\sigma_i) - \sum \frac {n_i}2 + \frac {\sum{\a_{ij}}}{2\pi} \\
&=& \sum \chi(\sigma_i) - E + V \\
&=& \chi(F) \qquad \qquad \qquad \qquad \qquad \qquad \hfill{  \Box}
\end{eqnarray*}
\bigskip

A face $\sigma$ is said to be spherical, euclidean or hyperbolic
according to whether $e(\sigma)$ is positive, zero or negative,
respectively.  When the angles are nonzero, one can show that $\sigma$
has a corresponding geometric structure with geodesic boundary edges
and an angle of $\a_i$ at corner $v_i$.  Thus for example if all
$\sigma_i$ are euclidean or hyperbolic with at least one $\sigma_i$
hyperbolic, and if the sum of angles at each vertex is at least
$2\pi$, then the surface $F$ is a hyperbolic surface.

\section{Proof of the main theorems}

Suppose $K(r)$ is an atoroidal Seifert fiber space.  By Lemma 6.2
there is an immersed sphere or a torus $\hat F$ in $K(r)$ which is in
essential position with respect to $P = P_1 \cup P_2 \cup P_3$.  Let
$G$ be the generalized graph on $\hat F$ as defined in Section 6.  Let
$\a_i$ and $\b_i$ be the angles of large corners and small corners of
$F\cap E(t_i)$, respectively.  We would like to show that if $\a_i,
\b_i$ can be chosen to satisfy certain conditions then $K(r)$ must be
toroidal.

As before, denote by $\bar \a_i = \pi - \a_i$ and $\bar \b_i = \pi -
\b_i$.  Recall that $\bar p_i$ is the mod $q_i$ inverse of $-p_i$, as
in Definition 2.2.  Let $\mu$ be the  meridional slope of $K$.  

%8.1
\begin{thm} Suppose $K = K(p_1/q_1, p_2/q_2, p_3/q_3)$ is a hyperbolic
  Montesinos knot with $q_i \geq 2$.  Then $K(r)$ is not an atoroidal
  Seifert fibered manifold if there are angles $\pi \geq \bar \a_i > 0$ and
  $\pi > \bar \b_i > 0$ satisfying the following conditions.

  (1) $\bar \a_1 + \bar \a_2 + \bar \a_3 \leq 2 \pi$;

  (2) $\bar \b_1 + \bar \b_2 + \bar \b_3 \leq \pi$;

  (3) $\bar \a_i + q_i\bar \b_i \geq 2\pi$;

  (4) $\bar \a_i +  |\bar p_i| \bar \b_i \geq \pi$;

  (5) if $q_i = 2$ then $\bar \a_i + \bar \b_i > \pi$.
\end{thm}

\proof Assume to the contrary that $K(r)$ is an atoroidal Seifert
fibred manifold.  Let $\hat F$ and $F$ be as above, with $\a_i$ and
$\b_i$ the large corners and small corners of $F \cap E(t_i)$,
respectively.  Conditions (1) and (2) can be rewritten as $\sum 2\a_i
\geq 2\pi$ and $\sum \b_i \geq 2\pi$, which mean that the sum of
the angles around each vertex of $G$ is at least $2\pi$.  Conditions
(3) and (4) say that the sum of external angles of a face $\sigma$ of
type $(1,q)$ or $(2,2|\bar p_i|)$ is at least $2\pi$, so by definition
we have $e(\sigma) \leq 0$.  We want to show that $e(\sigma)\leq 0$
for all other faces as well.

By Lemmas 3.3 and 3.4 a disk face of type $(r,s)$ in $E(t_i)$
is of one of the following type.

\qquad (a) $r = 0$ and $s\geq 2q_i$; 

\qquad (b) $r$ is odd and $s\geq q_i$;

\qquad (c) $r \geq 2$ is even and $s \geq 2|\bar p_i|$.

Let $\sigma_j$ be a face of type $j=a$, $b$, or $c$ above in $E(t_i)$.
Then
\begin{eqnarray*}
& & e(\sigma_a) = 1 - \frac 1{2\pi}\, s \bar \b_i \leq \frac
2{2\pi}(\pi - q_i \bar \b) \leq \frac{\bar \a_i - \pi}{\pi} \leq 0 \\ 
& & e(\sigma_b) = 1 - \frac 1{2\pi} (r \bar \a_i + s \bar \b_i) 
\leq 1 - \frac 1{2\pi} ( \bar \a_i + q_i \bar \b_i) \leq 0 \\
& & e(\sigma_c) = 1 - \frac 1{2\pi} (r \bar \a_i + s \bar \b_i) 
\leq  1 - \frac 1{2\pi} (2 \bar \a_i + 2|\bar p_i| \bar \b_i) \leq 0
\end{eqnarray*}

By definition all non-disk faces $\sigma$ of $G$ have $e(\sigma) \leq
0$.  Thus all faces $\sigma$ of $G$ have $e(\sigma) \leq 0$.  Since
$\chi(\hat F) = \sum e(\sigma) \leq 0$ by Lemma 7.2, the surface $\hat
F$ cannot be a sphere, so it must be a torus, and $e(\sigma) = 0$ for
all faces $\sigma$.  If there is a non-disk face then the outermost
one has some corner on it and hence has $e(\sigma)<0$, which is a
contradiction.  Therefore all $\sigma$ are disk faces with $e(\sigma)
= 0$.  This implies that all the above inequalities must be
equalities.  Since $\bar \a_i, \bar \b_i > 0$,  checking the above
inequalities gives

(a) If $\sigma_a$ exists in $E(t_i)$ then $\bar \alpha_i = \pi$, and
$s= 2q_i$, so $\sigma_a$ is a $(0, 2q_i)$ face.

(b) If $\sigma_b$ exists in $E(t_i)$ then $r=1$ and $s = q_i$, so
$\sigma_b$ is a $(1, q_i)$ face.

(c) If $\sigma_c$ exists in $E(t_i)$ then $r=2$ and $s = 2|\bar p_i|$,
so $\sigma_c$ is a $(2, 2|\bar p_i|)$ face.

These implies that $F$ satisfies condition (1) in Definition 5.1.
Condition (1) above implies that $\bar \a_i < \pi$ for some $i$, hence
(a) above implies that there is no face of type $(0, 2q_i)$ in $F \cap
E(t_i)$ for some $i$, so condition (2) in Definition 5.1 holds.  When
$q_i = 2$ we have $\bar p_i = 1$, and by condition (5) above we have
$\bar \a_i + \bar \b_i > \pi$, so by the calculation above we would
have $e(\sigma_c) < 0$, which is a contradiction.  Therefore there is
no disk of type $(2, 2) = (2, 2|\bar p_i|)$, which gives (3) in
Definition 5.1.  Similarly if $F \cap E(t_i)$ has a disk of type
$(0,2q_i)$ then $e(\sigma_a)=0$ gives $\bar \a_i = \pi$, and again we
have $\bar \a_i + |\bar p_i| \bar \b_i > \pi$, hence there is no disk
of type $(2, 2|p_i|)$, which verifies condition (4) of Definition 5.1.
Therefore $F$ is an elementary surface.  It now follows from
Proposition 5.6 that $K(r)$ is toroidal, a contradiction.  \qed

\bigskip
\noindent
{\bf Proof of Theorem 1.1.}  
Let $K = K(p_1/q_1, p_2/q_2, p_3/q_3)$, with $2 \leq q_1 \leq q_2 \leq
q_3$.  Since $K$ is a knot, at most one $q_i$ is even.  Hence the
condition $\sum \frac 1{q_i-1} \leq 1$ implies that either $q_i \geq
4$ for all $i$, or $q_1 = 3$ and $q_2, q_3 \geq 5$, or $q_1 = 3$, $q_2
= 4$ and $q_3 \geq 7$.

If $q_i \geq 4$, let $\bar \a_1 = \bar \a_2 = \bar \a_3 =
\frac{2\pi}{3}$, and $\bar \b_1 = \bar \b_2 = \bar \b_3 = \frac
{\pi}{3}$.

If $q_1 = 3$ and $q_2, q_3 \geq 5$, let $(\bar \a_1, \bar \a_2, \bar
\a_3) = (\frac {\pi}{2}, \frac {3\pi}{4}, \frac {3\pi}{4})$, and
$(\bar \b_1, \bar \b_2, \bar \b_3) = (\frac {\pi}{2}, \frac {\pi}{4},
\frac {\pi}{4})$.

If $q_1 = 3$, $q_2=4$ and $q_3 \geq 7$, let $(\bar \a_1, \bar \a_2,
\bar \a_3) = (\frac {\pi}{2}, \frac {2\pi}{3}, \frac {5\pi}{6})$, and
$(\b_1, \b_2, \b_3) = (\frac {\pi}{2}, \frac {\pi}{3}, \frac
{\pi}{6})$.

One can easily check that the conditions (1)--(3) of Theorem 8.1 are
satisfied in each of the above cases.  Note that in all cases we have
$\bar \a_i + \bar \b_i = \pi$, so condition (4) holds because $|\bar
p_i| \geq 1$.  Conditions (5) holds trivially since $q_i>2$.  Theorem
1.1 now follows from Theorem 8.1.  \qed

\bigskip
 
By Theorem 1.1 if a Montesinos knot $K$ of length 3 admits atoroidal
Seifert fibered surgery then $K = K(\frac{p_1}{q_1}, \frac{p_2}{q_2},
\frac{p_3}{q_3})$, such that either $q_1 = 2$, or $(q_1, q_2) =
(3,3)$, or $(q_1, q_2, q_3) = (3,4,5)$.  The following theorem gives
further restrictions on $p_i$.  In this theorem it is not assumed that
$q_2 \leq q_3$, so $q_2$ in (4) below may be larger than $q_3$.  Two
knots $K, K'$ are equivalent if $K$ is isotopic to $K'$ or its mirror
image.

%8.2
\begin{thm} Suppose a Montesinos knot $K$ of length 3 admits an
  atoroidal Seifert fibered surgery.  Then $K$ is equivalent to one of
  the following knots.

  (1) $K(1/3,\, \pm 1/4,\,  p_3/5)$ and $p_3 \equiv \pm 1$ mod $5$;

  (2) $K(1/3,\,  \pm 1/3,\,  p_3/q_3)$ and $|\bar p_3| \leq 2$;

  (3) $K(1/2,\,  2/5,\,  p_3/q_3)$, $q_3 = 5$ or $7$, and $|\bar p_3| > 1$;

  (4) $K(1/2,\,  1/q_2,\,  p_3/q_3)$, $q_2 \geq 5$ and $|\bar p_3| \leq 2$;

  (5) $K(1/2,\,  1/3,\,  p_3/q_3)$ and $|\bar p_3| \leq 6$.
\end{thm}

\proof By Theorem 1.1, we have $(q_1, q_2, q_3) = (3,4,5)$,
$(3,3,q_3)$ or $(2,q_2, q_3)$.  It is easy to see that if $K$ is not
equivalent one of those listed in the theorem then it is listed in one
of the following five cases.  In each case we will list the angles
$\bar \a_i, \bar \b_i$.  These satisfy $\sum \bar \a_i = 2\pi$, $\sum
\bar \b_i = \pi$, and $\a_1 = \pi$ when $q_1 = 2$, so conditions (1),
(2) and (5) in Theorem 8.1 hold.  We leave it to the reader to check
that they satisfy conditions (3) $\bar \a_i + q_i \bar \b_i \geq 2\pi$
and (4) $\bar \a_i + |\bar p_i| \bar \b_i \geq \pi$ in Theorem 8.1.

\medskip
{\sc Case 1.}  {\em $(q_1, q_2, q_3) = (3,4,5)$, and $|\bar p_3| =
  2$}. 

Let $(\bar \a_1, \bar \a_2, \bar \a_3) = (\pi, \frac{2\pi}{3},
\frac{\pi}{3})$, and $(\bar \b_1, \bar \b_2, \bar \b_3) = (\frac
{\pi}{3}, \frac {\pi}{3}, \frac {\pi}{3})$.

\medskip
{\sc Case 2.}  {\em $q_1 = q_2 = 3$, and $|\bar p_3|\geq 3$.}  

We must have $q_3 \geq 7$ because $2|\bar p_3| \leq
q_3$ and $\bar p_3$ is coprime with $q_3$.  Let $(\bar \a_1, \bar
\a_2, \bar \a_3) = (\frac {7\pi}{8}, \frac {7\pi}{8}, \frac {\pi}{4})$
and $(\bar \b_1, \bar \b_2, \bar \b_3) = (\frac {3\pi}{8}, \frac
{3\pi}{8}, \frac {\pi}{4})$.  

\medskip
{\sc Case 3.} {\em $q_1 = 2$, and $|\bar p_i| > 1$ for $i=2,3$.}

Note that $q_2, q_3$ must be odd since $K$ is a knot.  By definition
of $|\bar p_i|$ we have $q_2, q_3 \geq 5$.  If $q_2, q_3
\geq 7$, let $(\bar \a_1, \bar \a_2, \bar \a_3) = (\pi, \frac
{\pi}{2}, \frac{\pi}{2})$ and $(\bar \b_1, \bar \b_2, \bar \b_3) =
(\frac{\pi}{2}, \frac{\pi}{4}, \frac{\pi}{4})$.  If $q_2=5$ and $q_3
\geq 9$, let $(\bar \a_1, \bar \a_2, \bar \a_3) = (\pi, \frac
{\pi}{3}, \frac{2\pi}{3})$ and $(\bar \b_1, \bar \b_2, \bar \b_3) =
(\frac{\pi}{2}, \frac{\pi}{3}, \frac{\pi}{6})$.  Both cases can be
excluded by Theorem 8.1.  Therefore up to equivalence we have $q_2 =
5$ and $q_3 = 5$ or $7$, as in (3).

\medskip
{\sc Case 4.} {\em $q_1 = 2$, $|\bar p_2| = 1$, $q_2 \geq 5$, $|\bar
  p_3| \geq 3$.}

By convention we have $2|p_i|\leq q_i$ for $i=1,2$, hence $p_2 = \pm
1$.  The condition $|\bar p_3| \geq 3$ implies that $q_3 \geq 7$.  Let
$(\bar \a_1, \bar \a_2, \bar \a_3) = (\pi, \frac {3\pi}{4},
\frac{\pi}{4})$ and $(\bar \b_1, \bar \b_2, \bar \b_3) = (\frac
{\pi}{2}, \frac{\pi}{4}, \frac{\pi}{4})$.

\medskip
{\sc Case 5.} {\em $q_1 = 2$, $q_2 = 3$, $|\bar p_3| \geq 7$.}

We have $q_3 \geq 15$.  Let $(\bar \a_1, \bar \a_2, \bar \a_3) = (\pi,
\frac{7\pi}{8}, \frac{\pi}{8})$ and $(\bar \b_1, \bar \b_2, \bar \b_3)
= (\frac{\pi}{2}, \frac{3\pi}{8}, \frac{\pi}{8})$.

\medskip
In all cases above, we see that $\bar \a_i, \bar \b_i$ satisfy all
conditions of Theorem 8.1, therefore by that theorem $K(r)$ is not an
atoroidal Seifert fibred manifold, a contraction.  \qed

%8.3
\begin{rem} {\rm The first two tangles of the knots in the above
    theorem are very simple.  The small values of $|\bar p_3|$ implies
    that the third tangles of those knots are also relatively simple.
    For example if $\bar p_3 = 2$ and if we write $p_3 = n q_3 + m$
    with $2|m| < q_3$, then $p_3 \bar p_3 \equiv m\bar p_3\equiv - 1$
    mod $q_3$ implies that $q_3 = \pm (m \bar p_3 + 1) = \pm 2m \pm
    1$, so we have a partial fraction decomposition $p_3/q_3 = n \pm
    1/(2 \pm 1/m)$.  The results above will be used in [Wu4] to
    reduce the classification problem of Seifert fibered surgeries on
    Montesinos knots to that on a few specific families of knots. }
\end{rem}

\bigskip

\noindent
Department of Mathematics,  University of Iowa,  Iowa City, IA 52242
\\
Email: {\it wu@math.uiowa.edu}

\enddocument